\documentclass[12pt]{article}
\usepackage{mathrsfs}
\usepackage{amssymb,amsfonts,amsmath,amsthm,cite}
\usepackage{epsfig}
\newtheorem{thm}{Theorem}[section]
\newtheorem{cor}[thm]{Corollary}
\newtheorem{conj}[thm]{Conjecture}

\makeatletter \@addtoreset{equation}{section}

\def\qed{\hfill \rule{4pt}{7pt}}

\begin{document}

{\Large
\begin{center}
 {\large\bf A  Proof of  Moll's Minimum  Conjecture }
\end{center}
}

\begin{center}
William Y. C. Chen$^1$ and Ernest X. W. Xia$^2$ \\
Center for Combinatorics, LPMC-TJKLC\\
Nankai University\\
 Tianjin 300071, P. R. China

Email: $^1$chen@nankai.edu.cn, $^2$xxw@cfc.nankai.edu.cn

\end{center}


\noindent {\bf Abstract.} Let $d_i(m)$ denote the coefficients of
the Boros-Moll polynomials. Moll's minimum conjecture states that
the sequence $\{i(i+1)(d_i^2(m)-d_{i-1}(m)d_{i+1}(m))\}_{1\leq i
\leq m}$ attains its minimum with $i=m$. This conjecture is a
stronger  than the log-concavity conjecture proved by Kausers and
Paule. We give a proof of Moll's conjecture by utilizing the spiral
property of the sequence $\{d_i(m)\}_{0\leq i \leq m}$, and the
log-concavity of the sequence $\{i!d_i(m)\}_{0\leq i \leq m}$.



\section{Introduction}

The objective of this note is to give a proof of Moll's conjecture
on the minimum value of a sequence involving the coefficients of the
Boros-Moll polynomials which arise
 in the evaluation  of
the following quartic integral, see,
\cite{Alvarez2001,Am,George1999-2,George1999-3,George2001
 ,George2004,Moll2002}. It has been shown that
  for any $a>-1$ and any
nonnegative integer $m$,
\begin{equation*}
 \int_{0}^\infty
\frac{1}{(x^4+2ax^2+1)^{m+1}}dx
=\frac{\pi}{2^{m+3/2}(a+1)^{m+1/2}}P_m(a),
\end{equation*}
 where
\begin{align}\label{E2}
P_m(a)=2^{-2m}\sum_{k}2^k{2m-2k
 \choose m-k}{m+k \choose k}(a+1)^k.
\end{align}
Write $P_m(a)$ as
\begin{equation} \label{pma}
P_m(a)=\sum_{i=0}^md_i(m)a^i.
\end{equation}
The polynomials $P_m(a)$ are called the
 Boros-Moll polynomials.
By \eqref{pma},   $d_i(m)$ can be expressed as
\begin{align}\label{Defi}
d_i(m)=2^{-2m}\sum_{k=i}^m2^k{2m-2k \choose m-k}{m+k \choose
k}{k\choose i}.
\end{align}

From the above formula (\ref{Defi}) one sees that the coefficients
$d_i(m)$ are positive.
 Boros and Moll \cite{George1999-2, George1999-3}
 have proved that for $m\geq 2$ the sequence
$\{d_i(m)\}_{0 \leq i \leq m}$ is unimodal and the maximum entry
appears in the middle, that is,
\[d_0(m)<d_1(m)< \cdots<d_{\left[\frac{m}{2}\right]-1}(m)
<d_{\left[\frac{m}{2}\right]}(m)>
d_{\left[\frac{m}{2}\right]+1}(m)>\cdots >d_m(m).\]
  Moll \cite{Moll2002}  conjectured that the sequence
  $\{d_i(m)\}_{0 \leq i \leq m}$ is
log-concave for $m\geq 2$. Kauers and Paule \cite{Kauers2006} have
proved this conjecture by using a  computer algebra approach. Chen
and Xia \cite{Chen2008}
 have shown
that the sequence
 $\{d_i(m)\}_{0 \leq i \leq m}$ satisfies the  strongly
  ratio monotone
 property which implies the log-concavity and
 the spiral property.
  Chen and Gu \cite{Chen-Gu-2008}
  have proved that
 the sequence $\{d_i(m)\}_{0 \leq i \leq m}$
 satisfies
  the reverse ultra log-concavity. They have also proved that the
  sequence $\{i!d_i(m)\}_{0\leq i\leq m}$ is log-concave.

In fact, Moll \cite{Moll2007, Moll2008} proposed a stronger
conjecture than the log-concavity conjecture. He formulated his
conjecture in terms of the numbers $b_i(m)$  as defined by
\begin{align}\label{DE-B}
b_{i}(m)= \sum_{k=i}^m2^k{2m-2k \choose m-k} {m+k \choose
k}{k\choose i}.
\end{align}
Clearly, $b_i(m)=2^{2m}d_i(m)$ and the log-concavity of $d_i(m)$ is
equivalent to that of $b_i(m)$.

\begin{conj}\label{conj-1}
    Given $m\geq 2$, for $1\leq i \leq m$,
\[
(m+i)(m+1-i)b_{i-1}^2(m)+i(i+1) b_i^2(m)-i(2m+1)b_{i-1}(m)b_{i}(m),
\]
attains its minimum at $i=m$ with $2^{2m}m(m+1){2m \choose m}^2$.
\end{conj}

We will give a proof of the above conjecture by using the spiral
property of $\{d_i(m)\}_{0\leq i\leq m}$ and the log-concavity of
$\{i!d_i(m)\}_{0\leq i \leq m}$.

\section{Proof of  Moll's Minimum Conjecture}

As pointed out by Moll \cite{Moll2007}, his conjecture implies that
$\{d_i(m)\}_{0\leq i\leq m}$ is log-concave for $m\geq 2$. To see
this, we may employ a recurrence relation to reformulate his
conjecture by using the three terms $d_{i-1}(m)$, $d_i(m)$ and
$d_{i+1}(m)$. Recall that
  Kauers and Paule
\cite{Kauers2006} and Moll \cite{Moll2007} have independently
derived the following recurrence relation for $1\leq i \leq m$,
\begin{align}\label{recu4}
i(i-1)d_i(m)=(i-1)(2m+1)d_{i-1}(m)-(m+2-i)(m+i-1)d_{i-2}(m).
\end{align}
Note that we have adopted the convention that $d_i(m)=0$ for $i<0$
or $i>m$.
 From  \eqref{recu4} and the relation $d_i(m)=2^{-2m}b_i(m)$, it follows that
 \begin{align*}
(m+i)(m+1-i)b_{i-1}^2(m)&+i(i+1) b_i^2(m)
-i(2m+1)b_{i-1}(m)b_{i}(m)\\[6pt]
&=i(i+1)\left(b_{i}^2(m)-b_{i+1}(m)b_{i-1}(m)\right).
 \end{align*}
Thus, Moll's conjecture can be restated as follows.

\begin{thm}\label{Thm}
Given $m\geq 2$, for  $1 \leq i \leq m$,
  $i(i+1)\left(d_{i}^2(m)
 -d_{i+1}(m)d_{i-1}(m)\right)$ attains its
minimum at $i=m$ with
 $2^{-2m}m(m+1){2m \choose m}^2$.
\end{thm}

Chen and Xia
 \cite{Chen2008} have shown that the Boros-Moll polynomials satisfy
 the ratio monotone property which implies the log-concavity and the
 spiral property.

\begin{thm}
\label{Theo} Let $m \geq 2$ be an integer. The sequence
$\{d_i(m)\}_{0\leq i \leq m}$ is strictly
 ratio monotone, that is,
\begin{align*}
&\frac{d_m(m)}{d_0(m)}<\frac{d_{m-1}(m)}{d_1(m)}< \cdots
<\frac{d_{m-i}(m)}{d_i(m)}<\frac{d_{m-i-1}(m)}{d_{i+1}(m)}<\cdots
<\frac{d_{m-\left[\frac{m-1}{2}\right]}(m)}
{d_{\left[\frac{m-1}{2}\right]}(m)}<1,
\\[6pt]
 &\frac{d_0(m)}{d_{m-1}(m)}
 <\frac{d_1(m)}{d_{m-2}(m)}
<\cdots<\frac{d_{i-1}(m)}{d_{m-i}(m)}
<\frac{d_i(m)}{d_{m-i-1}(m)}<\cdots
<\frac{d_{\left[\frac{m}{2}\right]-1}(m)}
{d_{m-\left[\frac{m}{2}\right]}(m)}<1.
\end{align*}
\end{thm}

As a consequence of Theorem \ref{Theo}, the spiral property of
$\{d_i(m)\}_{0\leq i \leq m}$ can be stated as follows.

\begin{cor}{\rm (Chen and Xia
 \cite{Chen2008})}
\label{improvement-1} For $m \geq 2$, the sequence
$\{d_i(m)\}_{0\leq i \leq m}$ is  spiral, that is,
\begin{align}\label{spiral}
d_m(m)<d_0(m)<d_{m-1}(m)<d_{1}(m) <d_{m-2}(m)< \cdots
<d_{\left[\frac{m}{2}\right]}(m).
\end{align}
\end{cor}

Chen and Gu \cite{Chen-Gu-2008} have shown that
 $\{i!d_{i}(m)\}_{0 \leq i \leq
m}$ is log-concave. This property can be recast in the following
form.

\begin{thm}
\label{Log} For  $m \geq 2$ and $1 \leq i \leq m-1$,
\begin{align}\label{log-1}
id_i^2(m)>(i+1)d_{i+1}(m)d_{i-1}(m).
\end{align}
\end{thm}

We are now ready to present a proof of Theorem \ref{Thm}.

 \noindent
{\it Proof.} First,  it follows from \eqref{Defi} that
\begin{align}\label{step1}
m(m+1)d_m^2(m)=2^{-2m}m(m+1){2m \choose m}^2.
\end{align}
We now proceed to show that for $1\leq i \leq m-1$,
\begin{align}\label{step2}
i(i+1)\left(d_{i}^2(m)
 -d_{i+1}(m)d_{i-1}(m)\right)>m(m+1)d_m^2(m).
\end{align}
We first consider the case
 $1 \leq i \leq m-2$.
By  \eqref{log-1},  we find  that
\begin{align}\label{relatiom-1}
i(i+1)\left(d_i^2(m)-d_{i+1}(m)d_{i-1}(m)\right)
>i(i+1)d_i^2(m)-i^2d_i^2(m)
=id_{i}^2(m).
\end{align}
 Using the spiral property \eqref{spiral}, we see that
  for  $1 \leq i \leq m-2$,
\begin{align}\label{relatiom-2}
id_i^2(m)\geq d_1^2(m)>d_{m-1}^2(m).
\end{align}
Combining  \eqref{relatiom-1} and  \eqref{relatiom-2}, we get
\begin{align}\label{I-3}
i(i+1)\left(d_i^2(m)-d_{i+1}(m)d_{i-1}(m)\right)
>d_{m-1}^2(m).
\end{align}
On the other hand, by direct computation we may deduce from
\eqref{Defi}  that
\begin{align}\label{relatiom-3}
d_{m-1}(m)=\frac{2m+1}{2}d_m(m).
\end{align}
By \eqref{I-3} and \eqref{relatiom-3}, we have for  $1 \leq i \leq
m-2$,
\begin{align}
i(i+1)&\left(d_i^2(m)-d_{i+1}(m)d_{i-1}(m)\right)\nonumber\\[6pt]
&>\left(\frac{2m+1}{2}\right)^2d_m^2(m)
>m(m+1)d_m^2(m),
\label{relatiom-4}
\end{align} and hence  \eqref{step2} is true for  $1\leq i \leq m-2$.
 It remains to consider the case $i=m-1$.
Again, by \eqref{Defi} we find that
\begin{align}
d_{m-1}(m)&=2^{-m-1}
(2m+1){2m\choose m},\label{I-1}\\[6pt]
d_{m-2}(m)&=2^{-m-2} \frac{(m-1) (4m^2+2m+1)}{2m-1}{2m\choose m}.
 \label{I-2}
\end{align}
From  \eqref{step1}, \eqref{I-1} and \eqref{I-2},
 we deduce  that
\begin{align}\label{relatiom-5}
m(m-1)&\left(d_{m-1}^2(m)-d_m(m)d_{m-2}(m)\right)\nonumber
\\[6pt]
=&m(m-1)2^{-2m}{2m\choose m}^2
\left(\frac{(2m+1)^2}{4}-\frac{(m-1)(4m^2+2m+1)}{4(2m-1)}\right)
\nonumber
\\[6pt]
=&\frac{m(4m^2+6m-1)}{4(2m-1)}m(m-1)2^{-2m}
{2m\choose m}^2\nonumber\\[6pt]
>&m(m+1)2^{-2m}{2m\choose m}^2
=m(m+1)d_m^2(m).
\end{align}
Thus \eqref{step2}
 holds for $i=m-1$, and so it holds for $1\leq i \leq m-1$.
 This completes the proof.  \qed

To conclude, we propose the
 following ratio monotonicity  conjecture.
  If it
were true, it would imply that the sequence
$\{i(i+1)(d_i^2(m)-d_{i+1}(m)d_{i-1}(m))\}_{1\leq i \leq m}$ is both
spiral and log-concave for $m\geq 2$.

\begin{conj}
The sequence $\{i(i+1)\left(d_i^2(m)
-d_{i+1}(m)d_{i-1}(m)\right)\}_{1\leq i \leq m}$ is strongly ratio
monotone.
\end{conj}

 For example, for $m=8$, we have
 \allowdisplaybreaks
\begin{align*}
P_{8}(a)=&\frac{4023459}{32768}+\frac{3283533}{4096}a+
\frac{9804465}{4096}a^2+\frac{8625375}{2048}a^3
+\frac{9695565}{2048}a^4\\[6pt]
&+\frac{1772199}{512}a^5+\frac{819819}{512}a^6+\frac{109395}{256}a^7
+\frac{6435}{128}a^8.
\end{align*}
Let $c_i=i(i+1)\left(d_i^2(8)-d_{i+1}(8)d_{i-1}(8)\right)$ for
$1\leq i \leq 8$. One can verify that
\begin{align*}
\frac{c_8}{c_1}<\frac{c_7}{c_2}
 <\frac{c_6}{c_3} <\frac{c_5}{c_4}
<1\qquad \mbox{and}  \qquad \frac{c_1}{c_7}<\frac{c_2}{c_6}
<\frac{c_3}{c_5}<1 .
\end{align*}

\vspace{0.5cm}
 \noindent{\bf Acknowledgments.}  This work was supported by
  the 973
Project, the PCSIRT Project of the Ministry of Education, the
Ministry of Science and Technology, and the National Science
Foundation of China.


\begin{thebibliography}{99}
\bibitem{Alvarez2001}
J. Alvarez, M. Amadis, G. Boros, D. Karp, V.H. Moll and L. Rosales,
An extension of a criterion for unimodality, Electron. J. Combin.
8(1) (2001) R30.


\bibitem{Am}
T. Amdeberhan and V.H. Moll, A formula for a quartic integral: a
survey of old proofs and some new ones,
 Ramanujan J. to appear.

\bibitem{George1999-2}
G. Boros and V.H. Moll, A sequence of unimodal polynomials,  J.
Math. Anal. Appl.  237 (1999) 272--285.

\bibitem{George1999-3}
G. Boros and V.H. Moll, A criterion for unimodality,   Electron. J.
Combin.  6 (1999) \# R3.

\bibitem{George2001}
G. Boros and V.H. Moll, The double square root, Jacobi polynomials
and Ramanujan's Master Theorem, J. Comput. Appl. Math.  130 (2001)
337--344.

\bibitem{George2004}
G. Boros and V.H. Moll, Irresistible Integrals, Cambridge University
Press, Cambridge,  2004.

\bibitem{Chen-Gu-2008}
  W.Y.C. Chen and C.C.Y. Gu,  The reverse ultra log-concavity
 of the Boros-Moll polynomials,
    submitted,  arXiv:math.CO/0809.0127.

\bibitem{Chen2008}
W.Y.C. Chen and E.X.W. Xia, The ratio monotonicity of Boros-Moll
polynomials, Math. Comput., to appear.




\bibitem{Kauers2006}
M. Kausers and P. Paule, A computer proof of Moll's log-concavity
conjecture,   Proc. Amer. Math. Soc.  135(12) (2007)
 3847--3856.

 \bibitem{Moll2008}
D.V. Manna and V.H. Moll, A remarkable sequence of integers,
 arXiv: 0812.3374.

\bibitem{Moll2002}
V.H. Moll,  The evaluation of integrals:
 A personal story,  Notices
Amer. Math. Soc.  49(3) (2002)  311--317.

\bibitem{Moll2007}
V.H. Moll, Combinatorial sequences arising from a rational integral,
Online J. Anal. Combin.  2 (2007),
 \# 4.





\end{thebibliography}
\end{document}